\documentclass[11pt]{amsart}
\usepackage{pstricks,pstcol,pst-plot}

\definecolor{Cyan}{cmyk}{1,0,0,0}
\definecolor{GoldenRod}{cmyk}{0,0.10,0.84,0}
\definecolor{myblue}{rgb}{0.66,0.78,1.00}
\definecolor{Lavender}{cmyk}{0,0.48,0,0}
\definecolor{CornflowerBlue}{cmyk}{0.25,0.13,0,0}
\definecolor{Red}{cmyk}{0,1,1,0}
\definecolor{DarkBlue}{cmyk}{1,1,0,0.20}

\parskip=\smallskipamount
\numberwithin{equation}{section}

\newtheorem{theorem}{Theorem}[section]
\newtheorem{lemma}[theorem]{Lemma}
\newtheorem{corollary}[theorem]{Corollary}
\newtheorem{proposition}[theorem]{Proposition}

\theoremstyle{definition}
\newtheorem{definition}[theorem]{Definition}

\newtheorem{problem}[theorem]{Problem}
\newtheorem{remark}[theorem]{Remark}

\newcommand{\C}{{\mathbb C}}
\newcommand{\N}{{\mathbb N}}

\newcommand{\Z}{{\mathbb Z}}

\newcommand{\R}{{\mathbb R}}

\newcommand{\cA}{\mathcal{A}}
\newcommand{\cB}{\mathcal{B}}

\newcommand{\cC}{\mathcal{C}}
\newcommand{\cF}{\mathcal{F}}
\newcommand{\cH}{\mathcal{H}}

\newcommand{\cJ}{\mathcal{J}}
\newcommand{\cL}{\mathcal{L}}

\newcommand{\cO}{\mathcal{O}}

\newcommand{\cP}{\mathcal{P}}
\newcommand{\cU}{\mathcal{U}}

\def\di{\partial}
\def\dibar{\overline\partial}
\def\bs{\backslash}

\def\e{\epsilon}

\newcommand\wt{\widetilde}
\newcommand\hra{\hookrightarrow}
\newcommand\la{\langle}
\newcommand\ra{\rangle}

\begin{document}
\title{Stein compacts in Levi-flat hypersurfaces}

%
%
\author{Franc Forstneri\v c}
\address{Institute of Mathematics, Physics and Mechanics, 
University of Ljubljana, Jadranska 19, 1000 Ljubljana, Slovenia}
\email{franc.forstneric@fmf.uni-lj.si}
\thanks{The first author was supported by grants P1-0291
and J1-6173, Republic of Slovenia.}

%
%
\author{Christine Laurent-Thi\'ebaut}
\address{Institut Fourier, UMR 5582 CNRS/UJF, BP 74, 38402 St Martin
  d'H\`eres Cedex, France}
\email{Christine.Laurent@ujf-grenoble.fr}

%
%
\date{\today} 

\begin{abstract}
We explore connections between geometric properties of the 
Levi foliation of a Levi-flat hypersurface $M$ and holomorphic 
convexity of compact sets in $M$, or bounded in part by $M$. 
Applications include extendability of Cauchy-Riemann functions, 
solvability of the $\dibar_b$-equation, approximation 
of Cauchy-Riemann and holomorphic functions, and global regularity 
of the $\dibar$-Neumann operator.
\end{abstract}
\maketitle

\section{Introduction} 
A real hypersurface $M$ in an $n$-dimensional complex manifold is said to be 
{\em Levi-flat} if it is foliated by complex manifolds of dimension $n-1$; 
this {\em Levi foliation} is as smooth as $M$ according to Barrett and 
Forn\ae ss \cite{BaF}. A Levi-flat hypersurface locally partitions the 
complex manifold in two subsets which are both pseudoconvex along the 
hypersurface. Such hypersurfaces have recently received a lot of attention; 
see e.g.\ \cite{GT}, \cite{MAOH}, \cite{Ne}, \cite{Oh}, \cite{Siu2}.  

In this paper we describe a connection between the geometric properties 
of the Levi foliation and the complex analytic properties of certain 
compact sets in a Levi flat hypersurface. We indicate applications 
ranging from the approximation of Cauchy-Riemann and holomorphic functions, 
the exact solvability of the $\dibar_b$-equation (without shrinking the domain),
to the global regularity of the $\dibar$-Neumann operator on certain pseudoconvex 
domains in $\C^n$ containing a Levi-flat patch in the boundary.

In many analytic problems it is important to know that a certain compact set 
in a complex manifold admits a basis of open Stein neighborhoods; such a set will 
be called a {\em Stein compact} \cite{GR}, \cite{Ho}. Sometimes one needs Stein 
neighborhoods with certain additional analytic or topological properties.
A theorem of Siu \cite{Siu1} implies that 
each topologically closed leaf in a Levi-flat hypersurface which is exhausted  
by a strongly plurisubharmonic function admits a basis of open Stein neighborhoods 
in the ambient complex manifold (see also Col\c toiu \cite{Co} and Demailly \cite{De}).
Globally the situation is more complicated and far from well understood. 

One possible approach is to look for a holomorphic vector field transverse to 
$M$ on a given compact; its flow translates $M$ to a family of Levi-flat hypersurfaces, 
thus providing Stein neighborhoods of the original set. 
This approach was used by Bedford and Forn\ae ss \cite{BeF} to find precise results 
for complex curves in pseudoconvex boundaries, but so far it has been less successful 
in the case of Levi-flat hypersurfaces. Here we use a different
method. Our first  main  result is the following.

%
%
%
%
\begin{theorem} 
\label{t1.1}
Assume that $M$ is an orientable Levi-flat hypersurface of class $\cC^3$ 
in a complex manifold $X$. Let $\rho$ be a strongly plurisubharmonic $\cC^2$ 
function in an open set $U\subset X$ such that the set
$A = \{x\in U \cap M\colon \rho(x)\le 0\}$ is compact.
(Such $A$ will be called a {\em compact strongly pseudoconvex set in $M$.})
If the Levi foliation of $M$ is defined in a neighborhood 
of $A$ by a nowhere vanishing {\em closed} one-form of class $\cC^2$
then $A$ is a Stein compact. 
\end{theorem}

If we assume in addition that 
$d(\rho|_{M})\ne 0$ on $\{x\in M\cap U\colon \rho(x)=0\}=bA$ then $A$ 
is even {\em uniformly $\cH$-convex}
(Theorem \ref{t1.4} and Proposition \ref{p3.3}).

Theorem \ref{t1.1} is proved in \S 3.
First we show that the existence of a closed one-form defining
the Levi foliation of $M$ is equivalent to the existence of a $\cC^3$ 
defining function $v$ for $M$ such that $dd^c v=2i\di\dibar\, v$ 
and its first order derivatives vanish at every point of $M$
(Proposition \ref{p3.1}).
Such {\em asymptotically pluriharmonic} defining function is
used to find a Stein neighborhood basis of $A$ (Propositions \ref{p3.3} 
and \ref{p3.5}). With more regularity of $M$ we also obtain transverse 
holomorphic vector fields (Proposition \ref{p7.1}). 
The restriction to compacts is essential; for example, 
the Levi-flat hypersurface $M=\C\times\R \subset \C^2$ 
does not admit a basis of Stein neighborhoods in $\C^2$ \cite{AN}.

For the record we state here several sufficient conditions 
for the existence of a nonvanishing {\em closed} one-form 
defining a codimension one foliation; they are obtained 
from the classical theory of Reeb, Haefliger, 
Sullivan, Thurston and others (\S 4--\S 6).

\begin{proposition}
\label{p1.2}
A transversely orientable codimension one foliation $\cL$ 
of class $\cC^{r}$ $(r\ge 2)$ on a $\cC^r$ manifold $M$
is defined in a neighborhood of a certain compact set 
$A\subset M$ by a closed nowhere vanishing one-form of 
class $\cC^{r-1}$ provided that any one of the 
following conditions holds:
\begin{itemize} 
\item[(a)] $\mathcal{L}$ is a {\em simple foliation} in 
a neighborhood of $A$ (every point admits a local transversal 
to $\cL$ intersecting each leaf at most once).
\item[(b)] There is an open neighborhood $U\subset M$ 
of $A$ such that every leaf of the restricted foliation
$\cL|_U$ is topologically closed in $U$.
\item[(c)] $\mathcal{L}$ has no nontrivial one-sided holonomy
(this is the case when $M$ is real analytic, or when the leaves are
simply connected), and each element of $\pi_1(M)$ is of finite order. 
\item[(d)] $M$ is compact and there is a compact 
leaf $L\in\mathcal{L}$ with $H^1(L,\R)=0$. 
\item[(e)] $H^1(A,\mathcal{C}^{r-1}_{\mathcal{L}})=0$ where 
$\mathcal{C}^{r-1}_{\mathcal{L}}$ is the sheaf of real valued 
$\mathcal{C}^{r-1}$ functions on $M$ which are constant on the
leaves of $\mathcal{L}$.
\end{itemize}
\end{proposition}

The sufficient conditions in Theorem \ref{t1.1} and Proposition \ref{p1.2} 
seem fairly close to optimal. In \S 8 we show, using an example due to Bedford and 
Forn\ae ss \cite[p.\ 21]{BeF}, that the conclusion of Theorem \ref{t1.1}
(and of most analytic applications mentioned below)
fails in general in the presence of a leaf with nontrivial 
infinitesimal holonomy. In that example $M$ is a real analytic 
Levi-flat hypersurface in $\C^*\times \C$ whose intersection with a certain 
compact strongly pseudoconvex domain fails to be a Stein compact,
and both sides of $M$ are {\em worm domains} in the sense of 
Diederich and Forn\ae ss \cite{DF}. 

In \cite{SS2} Straube and Sucheston constructed closed one-forms defining 
certain Levi foliations on a compact hypersurface with boundary
$M\subset \C^2$, assuming that the leaves are closed and the foliation 
satisfies certain geometric conditions at $bM$.
(See also \cite[p.\ 256, Proposition]{SS1}.)
Using ideas of Sullivan \cite{Su} they characterized 
the existence of a smooth closed one-form defining the foliation 
in terms of foliation currents associated to a transverse flow 
\cite[p.\ 152, Proposition 2 and Theorem 3]{SS2}.

Condition (c) in Proposition \ref{p1.2} implies the following.

\begin{corollary}
\label{c1.3}
If $M$ is a simply connected, real analytic, Levi-flat hypersurface 
in a complex manifold $X$ then every compact strongly pseudoconvex 
subset $A \subset M$ is a Stein compact in $X$.
\end{corollary}

We now indicate some analytic applications of Theorem \ref{t1.1}.  Our  original motivation  
were problems on extendibility of Cauchy-Riemann (CR) functions and solvability 
of the $\dibar_b$-equation. Consider the following setup:
\begin{enumerate}
\item $X$ is a Stein manifold of dimension $n\ge 2$,
\item $D\subset \subset X$ is a smoothly bounded strongly pseudoconvex domain, 
\item $M$ is a Levi-flat hypersurface in $X$ intersecting $bD$ transversely,
\item $A=M\cap \overline D$, and
\item $\Omega$ is a connected component of $D\backslash M$ 
such that $b\Omega=A\cup \omega$, where $\omega$ is an open connected 
subset of $bD$. 
\end{enumerate}

Let $d$ be a distance function on $X$ induced by a smooth 
Riemannian metric. A compact set $A$ in $X$ is {\em uniformly $\cH$-convex} 
if there are open Stein domains $U_\nu\subset X$ $(\nu =1,2,\ldots)$ 
and a number $c\ge 1$ such that
\begin{equation}
\label{uniform-convex}
	\bigl\{x\in X\colon d(x,A)< \frac{1}{c\nu} \bigr\} \subset U_\nu \subset 
        \bigl\{x\in X\colon d(x,A)< \frac{c}{\nu} \bigr\}, \quad \nu \in\N. 
\end{equation}
The following result follows from Theorem \ref{t1.1} and Proposition \ref{p3.3}.

\begin{theorem}
\label{t1.4}
Let $X,M,D,A$ be as in (1)--(4), with $bD$ of class $\cC^2$. 
If one of the conditions in Theorem \ref{t1.1} or Proposition
\ref{p1.2} holds then the set $A$, and also the closure of 
every connected component of $D\backslash M$, 
is uniformly $\cH$-convex (and hence a Stein compact).  
\end{theorem}

Assuming the above setup (1)--(5) we consider the following analytic conditions;
(A2) and (A3) are relevant only if $n\geq 3$.
\begin{description} 
  \item[(A1)] Every continuous CR function on $\omega$ extends to a (unique) continuous 
  function on $E(\omega) = \overline{\Omega} \setminus A$ 
  which is holomorphic in $\Omega$. 
  
  \item[(A2)] For every smooth $(0,q)$-form $f$ on $\omega$ 
  $(1\leq q\leq n-3)$ satisfying $\dibar_b f=0$  there is a smooth 
  $(0,q-1)$-form $u$ on $\omega$ such that $\dibar_b u=f$.  
  
  \item[(A3)]  For every smooth $(0,n-2)$-form $f$ on $\omega$ for which 
  $\int_\omega f\wedge\varphi=0$ for every smooth, $\dibar$-closed $(n,1)$-form 
  $\varphi$ on $X\setminus A$ such that ${\rm supp}\,\varphi\cap{\rm supp}f$ is compact 
  there is a smooth $(0,n-3)$-form $u$ on $\omega$ satisfying $\dibar_b u=f$. 
\end{description} 
 
We recall the following result of the second author.

\begin{theorem} 
\label{Christine}
{\rm (\cite{CL3})}
If $\overline{\Omega}$ is a Stein compact then 
\[
	(A1), (A2) {\ \rm and\ } (A3) {\ \rm hold\ } \Longleftrightarrow 
	H^{n,q}(A)=0 {\ \rm for\ all\ } 1\leq q\leq n-1.
\]
\end{theorem}

Since all Dolbeault cohomology groups of a Stein compact vanish,
Theorems \ref{t1.1}, \ref{t1.4} and \ref{Christine} together imply

\begin{corollary}
\label{c1.6}
Assuming the setup (1)--(5), each of the conditions in 
Theorem \ref{t1.1} and Proposition \ref{p1.2}
implies that {\rm (A1), (A2), (A3)} hold. 
\end{corollary}

Results on extendability of CR functions from caps $\omega\subset bD$
have been obtained in \cite{CL1}, \cite{CL2}, \cite{CLP}, \cite{Lu},
\cite{LT}. For solutions of the $\dibar_b$-equation on $\omega$
without estimates up to the boundary see \cite{CLL};
for $L^p$ and Sobolev estimates see M.-C.\ Shaw \cite{Shaw1}, \cite{Shaw2}.

Our next application depends on the results of 
Boas and Straube \cite{BS1}, \cite{BS2} on regularity of the $\dibar$-Neumann operator 
$N_q$ on smoothly bounded pseudoconvex domains $\Omega\subset \subset \C^n$.
(See also Straube \cite{St2005} and Straube and Sucheston \cite{SS1}, \cite{SS2}.
We thank E.\ Straube for pointing out this connection.)
Recall that $N_q$ is the inverse of the complex Laplacian 
$\Box = \dibar \, \dibar^*+\dibar^*\dibar$ acting on $(0,q)$-forms,
and the Bergman projection $P_{q}=I-\dibar^* N_{q+1} \dibar$ is 
the orthogonal projection from the space of $(0,q)$-forms with $L^2$-integrable 
coefficients onto the subspace of $\dibar$-closed forms.
(For the $\dibar$-Neumann theory see \cite{Chen-Shaw}, \cite{FK}, \cite{Kohn}.)
Let $W^s_{(0,q)}(\Omega)$ denote the Sobolev space of all
$(0,q)$-forms whose coefficients have partial derivatives 
of order $\le s$ in $L^2(\Omega)$.

%
%
%
%
\begin{theorem} 
\label{Neumann}
Let $\Omega\subset \subset \C^n$ be a smoothly bounded 
pseudoconvex domain. Assume that the (compact) set $A\subset b\Omega$, 
consisting of all points of infinite D'Angelo type \cite{DA}, satisfies the following
properties:
\begin{itemize}
\item[(i)]   $A$ is the closure of its interior in $b\Omega$, and 
\item[(ii)]  $A$ is a strongly pseudoconvex set in a smooth Levi-flat hypersurface $M\subset \C^n$.
(See Theorem \ref{t1.1}.)
\end{itemize}
If the Levi foliation of $M$ is defined in a neighborhood of $A$
by a smooth closed $1$-form (in particular, if one of the conditions
in Proposition \ref{p1.2} holds) then the $\dibar$-Neumann operator $N_q$ 
and the Bergman projection $P_q$ are continuous on $W^s_{(0,q)}(\Omega)$ 
when $0\le q\le n$ and $s\ge 0$.
\end{theorem}

Theorem \ref{Neumann} is proved in \S 7. The conclusion may fail if $A$ 
is not the closure of its relative interior in $b\Omega$; a well known 
example is the complex annulus in the boundary of the Diederich-Forn\ae ss 
worm domain \cite{DF} (Barrett \cite{Ba3}, Christ \cite{Christ}). 
See also the example in \S 8 below.

%
%
%

We conclude this introduction by two theorems on holomorphic approximation
which are proved in \S 7. The first of them was proposed by G.\ Tomassini
(private communication). 
Note that a function on a Levi-flat hypersurface is a CR function 
if and only if its restriction to every leaf is holomorphic. Any CR function 
on a CR submanifold admits local holomorphic approximations \cite{BT}.

\begin{theorem}
\label{app1}
Let $A\subset M \subset X=\C^n$ be as in Theorem \ref{t1.1}. 
Assume that $M$ is of class $\cC^{r}$, $r\ge n+2$, and its Levi foliation 
is defined in a neighborhood of $A$ by a closed one-form of class $\cC^{r-1}$.
Let $k\in\{0,1,\ldots,r-1\}$ and $0<\alpha<1$.
Every CR function of class $\cC^{k,\alpha}$ in a neighborhood of $A$ in $M$ 
can be approximated in the $\cC^{k,\alpha}$ sense 
by functions holomorphic in a neighborhood of $A$ in $X$.
\end{theorem}

We believe that a version of Theorem \ref{app1} can also be proved 
for Levi-flat hypersurfaces in an arbitrary complex manifold 
since $A$ admits Stein neighborhoods; however, it appears that the smoothness 
hypothesis on $M$ should be  higher, and for simplicity  we restrict ourselves
to $X=\C^n$.

The final result depends on a theorem of Forn\ae ss and Nagel \cite{FN};
we thank E.\ Straube for calling this to our attention.

\begin{theorem}
\label{app2}
Let $\Omega \subset\subset \C^n$ be a smoothly bounded
pseudoconvex domain and $A$ a compact subset of $b\Omega$ 
satisfying the following: 
\begin{itemize}
\item[(i)]  $A$ is the closure of its interior in $b\Omega$, 
\item[(ii)] $b\Omega \bs A$ is strongly pseudoconvex, and
\item[(ii)] $A$ is a strongly pseudoconvex set in a 
smooth Levi-flat hypersurface $M\subset \C^n$ (see Theorem \ref{t1.1}).
\end{itemize}
If the Levi foliation of $M$ is defined in a neighborhood of $A$
by a smooth closed $1$-form (in particular, if one of the conditions
in Proposition \ref{p1.2} holds) then 
\begin{itemize}
\item[(1)] $\overline \Omega$ admits a basis of open Stein neighborhoods in $\C^n$, and
\item[(2)] every function which is holomorphic in $\Omega$ and continuous
on $\overline \Omega$ can be approximated uniformly on $\overline\Omega$ 
by functions holomorphic in a neighborhood of $\overline \Omega$.
\end{itemize}
\end{theorem}

%
%
%
%
\section{Preliminaries}
We denote by $j^r \eta$ the $r$-jet extension of a function 
or a differential form $\eta$ on a manifold $X$.  (In local coordinates
on $X$, $j^r \eta$ is the collection of all partial derivatives of order $\le r$ 
of the components of $\eta$.)  The notation $j^r \eta|_M=0$ will mean that the $r$-jet
$j^r\eta$ (with respect to all variables in the ambient manifold $X$)
vanishes at all points $p\in M$.

We recall some relevant notions concerning 
codimension one foliations and Levi-flat hypersurfaces. 
Our general references for the theory of foliations 
will be \cite{CN}, \cite{CC}, \cite{Go1}, and \cite{Ton}.

Let $M$ be a real manifold of dimension $m$ and class $\cC^r$, 
$r\in \{1,\ldots,+\infty,\omega\}$, where $\cC^\omega$ means 
real analytic. 
A foliation $\cL$ of codimension one and class $\cC^r$ on $M$ 
is given by a {\em foliation atlas} $\cU=\{(U_j,\phi_j)\colon j\in J\}$ 
where $\{U_j\}_{j\in J}$ is a covering of $M$
by open connected sets, for every index $j\in J$ the map 
$\phi_j=(\phi'_j, h_j) \colon U_j \to P_j= P'_j\times I_j \subset \R^{m-1}\times \R$
is a $\cC^r$ diffeomorphism, and the transition maps are of the form 
\[
	\theta_{ij}(u,v) = \phi_i\circ\phi^{-1}_j(u,v)= (a_{ij}(u,v),b_{ij}(v))
\]
where $b_{ij}$ is a diffeomorphism between a pair of intervals in $\R$.
Each leaf $L\in \cL$ intersects $U_j$ in at most countably many 
{\em plaques} $\{x\in U_j\colon h_j(x)=c \in I_j\}$. 
The collection $\{b_{ij}\}$ is called a {\em holonomy cocycle}, 
or a {\em Haefliger cocycle}, determining $\cL$ \cite{CC}, \cite{Ha}.
A foliation is {\em transversely orientable} 
(resp.\ {\em transversely real analytic}) if it admits a 
foliation atlas in which all diffeomorphism $b_{ij}$ 
preserve the orientation of $\R$ (resp.\ are real analytic). 
A continuous function $u\colon M\to \R$ is a {\em first integral} for $\cL$ 
if $u$ has no local extrema and is constant on every leaf of $\cL$;
$u \in \cC^1(M)$ is a {\em noncritical first integral} if 
in addition $du\ne 0$ on $M$. 

A closed loop $\gamma$ in a leaf $L\in \cL$
determines a germ of a diffeomorphism $\psi_\gamma$ 
(the {\em holonomy} of $\gamma$) on any local
transversal $\ell\subset M$ at a point $x_0 \in \gamma$, 
depending only on the homotopy class $[\gamma]\in \pi_1(L,x_0)$.
The induced map of $\pi_1(L,x_0)$ to the group of 
germs of diffeomorphisms of $(\ell,x_0)$ is called 
the (germinal) {\em holonomy homomorphism} of $L$.
$\cL$ has {\em trivial leaf holonomy} if $\psi_\gamma$ is the 
germ of the identity map for any loop $\gamma$ in a leaf of $\cL$;
this holds in particular if all leaves are simply connected.

We say that $\cL$ admits (nontrivial) {\em one-sided holonomy} if 
the holonomy map $\psi_\gamma$ of some leaf $L$, defined on a local transversal 
$\ell$ at some $p\in L$, equals the identity map 
on one side of $p$ in $\ell$, but not on the other side. 
A transversely real analytic foliation has no nontrivial 
one-sided holonomy.

Let $X$ be a complex manifold with the complex structure operator $J$. 
The operator $d^c=-J^* d$ is defined on functions 
by $\langle d^c v, \xi\rangle = -\langle dv, J\xi\rangle$
for $\xi\in TX$. In local holomorphic coordinates
$(z_1,\ldots,z_n)$, with $z_j=x_j+iy_j$ and 
$J(\frac{\di}{\di x_j})=\frac{\di}{\di y_j}$,  we have 
$d^c v =\sum_{j=1}^n -\frac{\di v}{\di y_j}\, dx_j + \frac{\di v}{\di {x_j}}\, dy_j$.
Then $d=\di+\dibar$, $d^c = i(\dibar -\di)$, and
$
	d d^c = 2i \di\dibar = 
	\sum_{j,k=1}^n \frac{\di^2}{\di z_j \di \bar z_k} dz_j\wedge d\bar z_k
$	
is the {\em Levi form operator}. A function  $v$ on $X$ 
is {\em pluriharmonic} if and only if $dd^c v=0$.

Let $M=\{v=0\} \subset X$ be a real hypersurface defined  
by a real function $v\in \cC^r(X)$ $(r\ge 2)$ with $dv\ne 0$ on $M$. 
The set $T^\C M= TM\cap J(TM) = TM\cap\, \ker{d^c} v$ 
is a real codimension one subbundle of $TM$.
The {\em Levi form} of $M$ is the quadratic map
$T^\C M\to \R$ given by 
\[
	T^\C M\ni \xi  \to \langle dd^c v, \xi \wedge J\xi \rangle =- \la d^c v,[\xi,J\xi]\ra. 
\]
(The last expression, which follows by Cartan's formula, 
assumes that $\xi$ is a vector field tangent to $T^\C M$.)
The hypersurface $M$ is {\em Levi-flat} if this form vanishes identically; 
this is equivalent to $dd^c v \wedge d^c v=0$ on $TM$
which is the integrability condition for $T^\C M =\ker(d^c v|_{TM})$. 
By Frobenius' integrability theorem the latter is equivalent to the 
existence of a codimension one foliation of $M$ by submanifolds
tangent to $T^\C M$, that is, by complex hypersurfaces. 
For Levi-flat hypersurfaces of low regularity 
see \cite{Ai} and  \cite{Sh}.

%
%
%
\section{Asymptotically pluriharmonic functions and Stein compacts}
In this section we prove Theorem \ref{t1.1}. We begin by preparatory results.

\begin{proposition}
\label{p3.1}
Let $M$ be an oriented, closed, Levi-flat hypersurface of class $\cC^r$ $(r\ge 2)$ 
in a complex manifold $X$. The following are equivalent:
\begin{itemize}
\item[(i)] $M$ admits a defining function $v\in \cC^r(X)$
with $j^{r-2}(dd^c v)|_M=0$.  

\item[(ii)] The Levi foliation of $M$ is defined by a closed
one-form of class $\cC^{r-1}$.
\end{itemize}
If $M$ is simply connected then (i) and (ii) are further equivalent to  
\begin{itemize}
\item[(iii)] There exists $f=u+iv \in \cC^r(X)$ such that 
$M=\{v=0\}$, $dv\ne 0$ on $M$, and $j^{r-1}(\dibar f)|_M=0$.
\end{itemize}
\end{proposition}

The equivalences in Proposition \ref{p3.1} also hold in the real analytic category
(see Remark \ref{real-analytic}). 
In that case (i) means that $M$ is defined by a function $v$ which is pluriharmonic 
in a neighborhood of $M$ in $X$, and the function $f=u+iv$ in (iii) is holomorphic in a 
neighborhood of $M$. 

A function $f$ satisfying $j^{r-1}(\dibar f)|_M=0$ is said to be 
{\em asymptotically holomorphic} of order $r-1$ on $M$, and a function $v$ 
satisfying  $j^{r-2}(dd^c v)|_M=0$ is {\em asymptotically pluriharmonic} 
of order $r-2$ on $M$. Since integrability of the subbundle $T^\C M$ 
is equivalent to $d^c v \wedge dd^c v|_{TM} =0$ for some
(and hence any) defining function $v$ for $M$ (\S 2),
a hypersurface $M=\{v=0\}$ with a defining function 
satisfying $dd^c v|_{TM} =0$ is necessarily Levi-flat.

\begin{proof} (i)$\Rightarrow$(ii):
Denote by $\iota \colon M\hookrightarrow X$ the inclusion map.
Let $\cL$ denote the Levi foliation of $M$, with the tangent 
bundle $T\cL = T^\C M$. Assume that $v\in \cC^r(X)$ satisfies $dv|_M\ne 0$ and 
$dd^c v|_M=0$. From  
$
	0= \iota^*(dd^c v) = d_M(\iota^* d^c v)
$
we see that the one-form $\eta = \iota^* (d^c v)$ 
is closed on $M$. Also, $\ker \eta = T^\C M = T\cL$
and hence $\eta$ defines the Levi foliation. 
Indeed,  for $\xi\in T_xM$ we have 
$\langle \eta, \xi \rangle = \langle d^c v, \xi \rangle = 
- \langle dv, J\xi\rangle$ 
which is zero if and only if $J\xi\in \ker dv_x=T_x M$, i.e., 
when $\xi\in T_x^\C M$. 

\smallskip
(ii)$\Rightarrow$(i):
Let $T\cL=\ker \eta$ where $\eta$ is a closed one-form of class $\cC^{r-1}$ 
on $M$. Locally near a point $p\in M$ we have $\eta=-du$ for some 
$\cC^r$ function $u$ which is unique up to an additive real
constant. Denoting by $\iota_L \colon L\hra M$  the inclusion of a 
leaf, we have $d(u|_L)=\iota_L^*(du)= \iota_L^* \eta=0$.
This shows that $u$ is constant on each leaf of $\cL$ and 
hence a Cauchy-Riemann (CR) function on $M$. Hence $u$ admits a
$\cC^r$ extension $f=u+iv$ to a neighborhood $U\subset X$ of 
$p$ such that $j^{r-1} \dibar f|_{M\cap U}=0$; such $f$ 
is unique up to a term which is flat to order $r$ along $M$
(\cite[p.\ 147]{Bo},  Theorem 2 and the remark following it).
Since $dd^c=2i\di\dibar$, the above implies 
$j^{r-2} dd^c(u+iv)|_{M\cap U}=0$. As $dd^c$ is a real operator, 
it follows that $j^{r-2} dd^c u|_{M\cap U}=0$ and 
$j^{r-2} dd^c v|_{M\cap U}=0$. Thus $v$ satisfies condition (i) on $M\cap U$. 
Since any two local solutions $v$ obtained in this way differ 
only by a term which is $r$-flat on $M$, we obtain 
a global solution $v$ by a smooth partition of unity.  
In the real analytic case the local pluriharmonic solutions 
are unique and no patching is needed.

If $M$ is simply connected and $T\cL=\ker \eta$ for some
one-form $\eta$ of class $\cC^{r-1}$ with $d\eta=0$, 
we have $\eta=-du$ for some $u\in \cC^r(M)$.
The same argument as above gives an asymptotically
holomorphic extension $f=u+iv$ into a neighborhood of $M$
in $X$ (holomorphic in the real analytic case),
thus showing the equivalence of (iii) with (i) and (ii). 
\end{proof}

\begin{remark}
\label{real-analytic}
The equivalence (i)$\Leftrightarrow$(ii) in Proposition \ref{p3.1}
was proved in the real analytic case by D.\ Barrett 
\cite[p.\ 461, Proposition 1]{Ba2}. For $r=2$ it also follows from 
results of Straube and Sucheston (combine the equivalence (i)$\Leftrightarrow$(iv) 
in \cite[p.\ 252, Theorem]{SS1} (in the exact case with $\epsilon=0$) 
and \cite[p.\ 152, Proposition 2]{SS2}). 

The equivalent conditions in Proposition \ref{p3.1}
hold in a small neighborhood of any point in $M$:
take a local $\cC^r$ first integral $u$ 
of $\cL$ and extend it to a $\dibar$-flat function 
$f=u+iv$ as above. Such $f$ gives a local 
asymptotically holomorphic flattening of $M$ which 
in general cannot be chosen holomorphic on any side 
of $M$ by \cite{BdB}. A real analytic Levi-flat hypersurface
need not admit a holomorphic flattening on large domains,
not even on contractible ones \cite{F2}. 
\end{remark}

Theorem \ref{t1.1} will follow from Proposition \ref{p3.1}
and the following result. 

\begin{proposition}
\label{p3.3}
Let  $M$ be a Levi-flat hypersurface of class $\cC^3$ 
in a complex manifold $X$ and $\rho$ a function of class $\mathcal{C}^2$ 
in an open set $U\subset X$ such that 
\begin{itemize}
\item[(i)] the set $A = \{x\in U \cap M\colon \rho(x)\le 0\}$ is compact,
\item[(ii)] $d(\rho|_M) \ne 0$ at every point of $bA=\{x\in M\cap U\colon \rho(x)=0\}$, and 
\item[(iii)] $dd^c \rho>0$ on $T^\C_x M$ at each point $x\in A$.
\end{itemize}
If $M\cap U$ admits a defining function $v\in \cC^3(U)$ 
satisfying $j^1(dd^c v)|_{A}=0$ then $A$ is uniformly 
$\cH$-convex, and hence a Stein compact.
\end{proposition}

\begin{proof}
Choose a smooth Riemannian metric on $X$ and denote by $|\xi|$ 
the associated norm of any tangent vector $\xi \in TX$. 
Since $dd^c\rho >0$ on $T^\C M$ at the points of $A$, 
the function $\wt \rho = \rho+ Cv^2$ is strongly plurisubharmonic 
in an open neighborhood of $A$ in $X$ provided that the constant
$C>0$ is chosen sufficiently large. Replacing $\rho$ by $\wt \rho$ 
and shrinking $U$ around $A$ we may therefore assume that 
the quantities $|dv|$ and $dd^c \rho >0$ are bounded, and also
bounded away from zero, on $U$. The distance of a point $x\in U$ to $M$ 
is comparable to $|v(x)|$.

Set $U^\pm  =\{x\in U \colon \pm \, v(x)\ge 0\}$. 
Choose a small $c>0$ such that 
$\{x\in U\cap M \colon  \rho(x) \le c\} \subset\subset  U$
and each value in $[0,c]$ is a regular value of $\rho|_{M\cap U}$.
Then $dv$ and $d\rho$ are $\R$-linearly independent 
at every point of $E=\{x\in M\cap U\colon \rho(x)=c\}$.
Choosing $U$ sufficiently thin around $M\cap U$ we 
insure that the set $D=\{x\in U\colon \rho(x)\le c\}$ has $\cC^2$
boundary intersecting $M$ transversely along $E$. 
For small $\epsilon>0$ and $x\in U$ set 
\begin{eqnarray*}
	v_\epsilon^\pm(x) &=& \pm\, v(x) + \epsilon (\rho(x)-c), \cr 
        \Gamma_\epsilon^\pm &=& \{x\in U \colon   v_\epsilon^\pm(x)=0\}, \cr
        \Omega_\epsilon &=& \{x\in U\colon v^+_\epsilon(x)<0,\ v^{-}_\epsilon(x)<0,\ 
        \rho(x)<\epsilon \}. 
\end{eqnarray*}
(See fig.\ \ref{Fig1}.)
For sufficiently small $\epsilon >0$ we have
$dv_\epsilon^\pm= \pm dv +\epsilon d\rho\ne 0$ on $U$;
hence $\Gamma_\epsilon^\pm$ are $\cC^2$ hypersurfaces 
satisfying $\Gamma_\epsilon^\pm \cap  D \subset U^\pm$ and
$\Gamma_\epsilon^+\cap \Gamma_\epsilon^- = E$.
As $\epsilon \to 0$, the domains $\Omega_\epsilon$ 
shrink down to $A$.

%
%
%
%
%
%

\begin{figure}[ht]
\psset{unit=0.7cm} 
\begin{pspicture}(-9,-4)(9,4)

\psline[linewidth=0.8pt](-9,0)(9,0)                                
\psline[linewidth=2.5pt,linecolor=Red](-5.2,0)(5.2,0)              
\psdots[dotsize=4pt,linecolor=Red](-5.2,0)(5.2,0)

\psarc[linestyle=dotted,linewidth=1pt](0,-12){14}{59}{121}         
\psarc[linewidth=1.4pt,linecolor=DarkBlue](0,-12){14}{65}{115}                  

\psarc[linestyle=dotted,linewidth=1pt](0,12){14}{-121}{-59}        
\psarc[linewidth=1.4pt,linecolor=DarkBlue](0,12){14}{-115}{-65}

\psarc[linestyle=dotted,linewidth=1pt](-2,0){4}{140}{220}
\psarc[linewidth=1.4pt,linecolor=DarkBlue](-2,0){4}{170}{190}      

\psarc[linestyle=dotted,linewidth=1pt](2,0){4}{-40}{40}
\psarc[linewidth=1.4pt,linecolor=DarkBlue](2,0){4}{-10}{10}        

\rput(0,1){$\Omega_\e$}

\rput(-7.8,-1.2){$M$}
\psline[linewidth=0.2pt]{->}(-7.8,-0.8)(-7.8,-0.07) 

\rput(-2.5,1){$A$}
\psline[linewidth=0.2pt]{->}(-2.5,0.7)(-2.5,0.07) 

\rput(3,2.8){$\Gamma_\e^+$}
\psline[linewidth=0.2pt]{->}(3,2.4)(3,1.8) 

\rput(3,-2.8){$\Gamma_\e^-$}
\psline[linewidth=0.2pt]{->}(3,-2.4)(3,-1.8)

\rput(-7.2,2){$\{\rho=\e\}$}
\psline[linewidth=0.2pt]{->}(-6.2,2)(-5.5,2)

\rput(7.2,2){$\{\rho=\e\}$}
\psline[linewidth=0.2pt]{<-}(5.5,2)(6.2,2)

\psdot[dotsize=4pt](-7.2,0)
\psdot[dotsize=4pt](7.2,0)

\rput(7.2,-1.2){$E$}
\psline[linewidth=0.2pt]{->}(7.2,-0.8)(7.2,-0.15)

\end{pspicture}
\caption{The domain $\Omega_\epsilon$} 
\label{Fig1}
\end{figure}
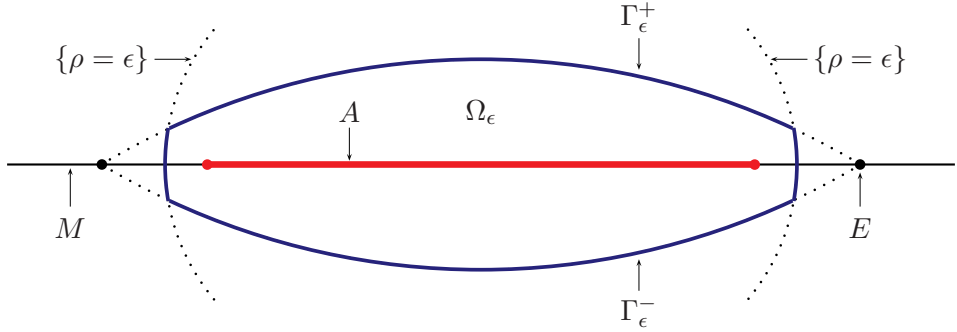

We claim that the functions $v_\epsilon^\pm$ are strongly plurisubharmonic 
on $\Omega_\epsilon$ provided that $\epsilon>0$ is sufficiently small. 
By the choice of $U$ there is $C_0>0$ such that 
$\langle dd^c\rho(x), \xi \wedge J\xi \,\rangle \ge C_0 |\xi|^2$ 
for $x\in U$ and $\xi\in T_x X$. Since $j^1(dd^c v)$ vanishes on $A$,  
we have $|dd^c v(x)| = o(|v(x)|+\epsilon)$. (The extra $\epsilon$ takes
care of the points near $\Omega_\e \cap M\bs A$.)
For $x\in  \Omega_\epsilon$ we also have 
\[
	|v(x)| \le \epsilon |\rho(x)-c| \le C_1\epsilon
\]
with $C_1=\sup \{|\rho(x)-c|\colon x\in D\}$;
hence $|dd^c v(x)| = o(\epsilon)$ for $x\in\Omega_\epsilon$ as $\epsilon\to 0$. 
This gives for $x\in \Omega_\epsilon$ and $\xi \in T_x X$:
\begin{eqnarray*}
	\langle dd^c v_\epsilon^\pm(x), \xi\wedge J\xi \,\rangle 
	&=& \epsilon \langle dd^c \rho(x),  \xi \wedge J\xi \,\rangle 
	    \pm  \langle dd^c v(x), \xi\wedge J\xi \,\rangle\\ 
        & & \ge \bigl( C_0\epsilon -o(\epsilon)\bigr) \cdotp |\xi|^2
\end{eqnarray*}
which is positive for sufficiently small $\epsilon>0$ and  
for $\xi\ne 0$. This establishes the claim.

Choose a smooth strongly increasing convex
function $h\colon (-\infty,0)\to \R$ with 
$\lim_{t\to 0} h(t)=+\infty$. The function 
\[ 
	\tau(x)= h(v_\epsilon^+(x))+h(v_\epsilon^-(x)) + h(\rho(x)-\epsilon), \quad x\in\Omega_\e 
\]
is then a strongly plurisubharmonic exhaustion function on 
$\Omega_\epsilon$, and hence $\Omega_\epsilon$ is Stein according to
\cite[p.\ 127, Theorem 5.2.10]{Ho}. The family $\{\Omega_\e\}$ 
is therefore a Stein neighborhood basis of $A$ in $X$ satisfying (\ref{uniform-convex}).
\end{proof}

\begin{remark}
\label{blowingup}
If $v\in\cC^3(U)$ is a defining function for $M\cap U$ which
satisfies $j^1(dd^c v)|_{M\cap U}=0$ then the function 
$\tau=-\log v +\rho$, defined on $U\cap \{v>0\}$ 
and tending to $+\infty$ along $M\cap U$, is strongly plurisubharmonic 
near $A$. Indeed, we have
\[
	dd^c \tau = -\frac{1}{v}\, dd^c v +
	\frac{1}{v^2}\, dv\wedge d^c v + dd^c \rho.
\]
From $j^1(dd^c v)|_M=0$ we infer  that the first term vanishes on $M$,
the second term is nonnegative since 
\[
	\la dv\wedge d^c v, \xi\wedge J\xi\ra = 
	\la dv,\xi\ra \la d^c v,J\xi\ra - \la dv,J\xi\ra \la d^c v,\xi\ra
	= \la dv,\xi\ra^2 + \la dv,J\xi\ra^2,
\]
and the last term $dd^c\rho$ is positive since $\rho$ is strongly plurisubharmonic.
The analogous observation holds for $-\log(-v) +\rho$
from the side $\{v<0\}$.
\end{remark}

\textit{Proof of Theorem \ref{t1.1}.}
Choose a small $c>0$ which is a regular value of the function $\rho|_M$
such that the set $A_c= \{x\in M\cap U \colon \rho(x)\le c\}$
is compact and the hypothesis concerning the Levi foliation of $M$ 
still holds on a neighborhood of $A_c$. By Proposition \ref{p3.1}
there exists a $\cC^3$ defining function $v$ for $M$ near $A_c$
satisfying $j^1(dd^c v)=0$ on $A_c$. Proposition \ref{p3.3} now implies 
that $A_c$ is uniformly $\cH$-convex and hence a Stein compact. 
Since there exist arbitrarily small numbers $c>0$ for which the above 
holds, the set $A=A_0$ is a Stein compact as well. 
This concludes the proof of Theorem\ref{t1.1}.
\smallskip

Theorem \ref{t1.4} is obtained in exactly the same way as
Theorem \ref{t1.1} by combining Proposition \ref{p3.1}
and the following result.

\begin{proposition}
\label{p3.5}
Let $D$ be a relatively compact, strongly pseudoconvex domain in a Stein manifold $X$, 
and let $M \subset X$ be a $\cC^3$-smooth Levi-flat hypersurface intersecting 
$bD$ transversely. Set $A=M\cap \overline D$. If $M$ admits a defining function 
$v\in\cC^3(X)$ satisfying $j^1(dd^c v)|_A=0$ then the set $A$, as well as the 
closure of any connected component of $D\bs M$, is uniformly $\cH$-convex,
and hence a Stein compact in $X$. 
\end{proposition}

\begin{proof} 
Let $D=\{\rho<0\}$ where $\rho$ is a strongly plurisubharmonic 
function in a neighborhood of $\overline D$ in $X$.
Assume first that $D\bs M=D^+\cup D^-$ 
consists of two connected components $D^\pm$.
We may assume that $v>0$ in $D^+$ locally near $M\cap D$. 
For sufficiently small $\epsilon>0$ the hypersurface 
$\Gamma_\epsilon^-$, constructed in the proof of Proposition  
\ref{p3.3}, intersects the hypersurface $\{\rho=\epsilon\}$ transversely. 
Let $\Omega_\epsilon$ denote the domain bounded by 
$\Gamma_\epsilon^-$ and by the part of $\{\rho=\epsilon\}$
on which $v_\epsilon^- <0$. (Thus $\Omega_\epsilon$ contains
$\overline{D^+}$ but not $\overline{D^-}$.) 
The argument in the proof of Proposition \ref{p3.3} 
shows that $\Gamma_\epsilon^-$ is strongly pseudoconvex
and hence $\Omega_\epsilon$ is Stein.
Thus $\overline {D^+}$ is uniformly $\cH$-convex. 
The analogous argument applies to any connected component 
of $M\cap D$, and also to $A$ itself.  
\end{proof}

%
%
%
%
%
\section{Simple foliations}
In this and the following two sections we show how
Proposition \ref{p1.2} follows from the classical theory of 
codimension one foliations. We provide 
sufficient details and references to make the paper accessible 
to a reader without a substantial background in foliation theory.

By $\cL$ we shall always denote a transversely orientable 
codimension one foliation on a manifold $M$. 
We begin by recalling some basic facts.

1. $\cL$ is defined by a closed one-form if and only if 
it admits a foliation atlas whose holonomy cocycle
consists of translations in $\R$ \cite[p.\ 28]{CC}, \cite{Ha}. 

2. A foliation defined by a closed one-form $\eta$ 
has trivial leaf holonomy \cite[p.\ 80]{CN}, \cite[Theorem 3.29]{Ton}.
Indeed, if $\gamma$ is a simple closed loop in a leaf 
$L$ then $\int_\gamma \eta=0$ (since $\eta|_L=0$);
choose $p\in\gamma$ and define $u(x)=\int_{p}^x \eta$ 
(the integral is independent of the path in a tubular 
neighborhood of $\gamma$);
then $du=\eta \ne 0$ and hence $u$ is a first integral for 
$\cL$ which is injective on a local transversal at $p$, 
so the germinal holonomy of $L$ along $\gamma$ is trivial. 
Conversely, a transversely orientable codimension one 
foliation $\cL$ without holonomy and of class 
$\cC^r$ $(r\in \{2,\ldots,+\infty\})$ is topologically, 
but in general not diffeomorphically, conjugate to a $\cC^r$ foliation 
defined by a closed nowhere vanishing one-form (Sacksteder \cite{Sa},  
\cite[\S 9.2 and p.\ 218]{CC}).

3. The {\em Godbillon-Vey class} of a foliation $\cL$
given by a one-form $\eta$ is 
$gv(\cL)=[\eta\wedge \alpha\wedge \beta]\in H^3(M,\R)$ where 
$\alpha$, $\beta$ are one-forms satisfying 
$d\eta=\alpha \wedge \eta$, $d\alpha =\beta \wedge \eta$ 
(\cite[p.\ 38]{CC}, \cite{Gh}, \cite{GV}, \cite[Theorem 2.3]{Ton}). 
Clearly $d\eta=0$ implies $gv(\cL)=0$.

%
%
\begin{definition}
\label{simplefol}
{\rm (\cite[p.\ 79]{Go1})}
A foliation $\cL$ on a manifold $M$ is {\em simple} 
if every point $p\in M$ is contained in a local
transversal $\ell \subset M$ to $\cL$ which is not intersected 
more than once by any leaf of $\cL$.
\end{definition}

We begin the proof of Proposition \ref{p1.2} by observing that 
conditions (a) and (b) are equivalent.
Indeed, choosing an open relative neighborhood 
$U\subset\subset M$ of $A$ with $\cC^1$ boundary, 
the first Betti number of $U$ is finite and hence the condition (b) 
(that all leaves of $\cL|_U$ are closed) implies that the foliation 
is simple \cite[p.\ 116, Lemma II.3.9]{Go1}. 
The converse is seen by observing that a nonclosed leaf must 
accumulate on some leaf $L$, and hence it will intersect a local
transversal to $L$ infinitely many times. 

The implication (d)$\Rightarrow$(a) 
follows from the stability theorems of Reeb and Thurston  
\cite[p.\ 78, Theorem 5]{CN}, \cite[\S 6.2]{CC}, \cite[II.3]{Go1}. 
The implication (c)$\Rightarrow$(a) will be shown in \S 5,
and (e) is discussed in \S 6.

The following result, which was proved for simply connected 
manifolds by Haefliger and Reeb \cite{HR}, shows the sufficiency 
of (a) in Proposition \ref{p1.2}.

\begin{theorem}
\label{t4.2}
A transversely orientable simple foliation of 
codimension one and class $\cC^r$ $(r\in \{2,\ldots, +\infty\})$ 
is defined in a neighborhood of any compact set 
by a closed nowhere vanishing one-form of class $\cC^{r-1}$. 
\end{theorem}

Globally such one-form need not exist even if $M=\R^2$ \cite{Wa},
and it cannot be chosen real analytic  even if
the foliation is real analytic \cite[p.\ 122]{Go0}. 

Having been unable to locate a precise reference for Theorem 
\ref{t4.2}, we include a proof following the methods of Haefliger and 
Reeb \cite{HR}.  We begin by recalling some results concerning 
the leaf space of a simple foliation.

\begin{proposition}
\label{p4.3}
{\rm (Haefliger and Reeb \cite{HR})}  
Let $\cL$ be a simple foliation of codimension one
and of class $\cC^r$ $(r\ge 1)$ on a connected manifold $M$. 
The space of leaves $Q=M/\cL$ carries the structure of
a connected, one dimensional,  not necessarily Hausdorff $\cC^r$ manifold  
such that the quotient projection $\pi\colon M\to Q$ is a $\cC^r$ submersion. 
The leaf space $Q$ is orientable if and only if $\cL$ is 
transversely orientable, and $Q$ is simply connected if and 
only if $M$ is such. Every function on $M$ which is constant 
on the leaves of $\cL$ is of the form $f\circ \pi$ for some 
function $f\colon Q\to \R$.
\end{proposition}

The $\cC^r$ structure on the leaf space $Q=M/\cL$ is 
uniquely determined by the requirement that the restriction 
of the quotient projection $\pi\colon M\to Q$ to any 
local $\cC^r$ transversal $\ell\subset M$ to the foliation $\cL$ 
is a $\cC^r$ diffeomorphism of $\ell$ onto $\psi(\ell) \subset Q$.

To prove Theorem \ref{t4.2} it suffices 
to find a nowhere vanishing one-form $\theta$ on the leaf
space $Q=M/\cL$; it's pull-back $\eta = \pi^* \theta$
is then a nowhere vanishing one-form on $M$ satisfying 
$\ker \eta = T\cL$. A non-Hausdorff manifold $Q$ need not 
admit any nonconstant $\cC^1$ functions 
or one-forms (Wazewski \cite{Wa}; \cite{HR}); 
however, one can do this on compact subsets 
as we now explain.

A point $q$ in a non-Hausdorff manifold $Q$ is said to be 
a {\em branch point} if there exists another
point $q'\in Q$ different from $q$ such that 
$q$ and $q'$ have no pair of disjoint neighborhoods; 
such $\{q,q'\}$ is called a {\em branch pair}.
(This relation is not transitive.) Branch points 
of $Q=M/\cL$ correspond to {\em separatrices}
of the foliation $\cL$. For a given $q\in Q$ 
there may exist infinitely many $q'\in Q$
such that $\{q,q'\}$ is a branch pair.  
The main difficulty is that
a germ of a smooth function at a branch point need 
not correspond to the germ of any global function on $Q$ 
(see an example in \cite{HR}).  A one dimensional 
non-Hausdorff manifold $Q$ of class $\cC^r$  is said to 
be {\em regular} if every germ of a $\cC^r$ function 
at any point $q\in Q$ is the germ at $q$ of a global 
$\cC^r$ function on $Q$ \cite[p.\ 116, Def.\ 2]{HR}. 
The following results are due to Haefliger and Reeb.

\begin{proposition} 
\label{p4.4}
{\rm (\cite[p.\ 123, Proposition]{HR})}
Let $\cL$ be a simple $\cC^r$ foliation of codimension one on a manifold $M$.  
For any open relatively compact domain $A\subset M$ the leaf space 
$Q_A=A/\cL_A$ of the restricted foliation $\cL_A$ is a regular,
not necessarily Hausdorff, one dimensional $\cC^r$ manifold.
\end{proposition}

\begin{proposition} 
\label{p4.5}
{\rm (\cite[p.\ 117, Proposition 1]{HR})}
A regular, simply connected, one dimensional manifold of class 
$\cC^r$ $(r\ge 1)$, not necessarily Hausdorff, admits a 
$\cC^r$ function with nowhere vanishing differential. 
\end{proposition}

Combining Propositions \ref{p4.3}, \ref{p4.4}, and \ref{p4.5}
one gets the following result which was proved for 
smooth foliations of $\R^2$ by Kamke \cite{Ka}.

\begin{corollary}
\label{c4.6}
Let $r\ge 1$. A simple foliation of codimension one and class $\cC^r$  
on a manifold $M$ admits a noncritical first integral of class $\cC^r$ 
on any relatively compact, simply connected open set in $M$.
\end{corollary}

Assume now that $(M,\cL)$ satisfies the hypotheses of Theorem \ref{t4.2}.
Choose a compact set $K$ in $M$. 
By \cite[Lemma IV.1.6]{Go1}  there is a relatively compact
domain $A\subset\subset  M$, with $K\subset A$, whose boundary
$bA$ is of class $\cC^r$ and {\em in general position} with respect 
to $\cL$, meaning that any local defining function 
for $\cL$ is a Morse function when restricted to $bA$, and 
distinct critical points of these functions belong to 
distinct leaves of $\cL$ \cite[p.\ 228, Def.\ IV.1.4.]{Go1}.

Note that a point $p\in bA$ is a critical point of a local defining 
function for $\cL$ (restricted to $bA$) precisely when the leaf
through $p$ meets $bA$ tangentially at $p$ (it is a separatrix).  
By our choice of $A$ all but finitely many leaves intersect 
$bA$ transversely, and hence the 
leaf space $Q_A=A/\cL_A$ has at most finitely many branch points. 
Furthermore, since distinct points of contact of the leaves with 
$bA$ belong to distinct leaves, $Q_A$ has at most double branch points.

To complete the proof of Theorem \ref{t4.2} it suffices
to show the following.

\begin{proposition}
\label{p4.7}
An oriented, not necessarily Hausdorff, one dimensional manifold
of class $\cC^r$ $(r\in \{1,2,\ldots,+\infty\})$ which is regular 
and has at most finitely many double branch points 
(and no multiple branch points) admits a nowhere 
vanishing differential one-form of class $\cC^{r-1}$.
\end{proposition}

\begin{proof}
Fix an orientation on $Q$.
Let $x_j,y_j \in Q$ $(j=1,\ldots, k)$ be all pairs of branch
points, i.e., $x_j$ and $y_j$ cannot be separated by open neighborhoods,
but any other pair of distinct points of $Q$ can be separated.
For every $j$ we choose an open neighborhood $U_j \subset Q$
of $x_j$ and an orientation preserving $\cC^r$ diffeomorphisms 
$h_j\colon U_j\to I=(-1,1)\subset \R$. By regularity
of $Q$ we can assume that $h_j$ extends to a $\cC^r$ function 
$h_j\colon Q\to \R$. The extended function necessarily has nonzero derivative 
at the point $y_j$ which forms a branch pair with $x_j$ \cite{HR}, 
and hence $\theta_j=dh_j$ is a nowhere vanishing 
one-form of class $\cC^{r-1}$ in an open connected neighborhood 
$V_j\subset Q$ of the pair $\{x_j,y_j\}$. We may assume that the 
closures $\overline V_j$ for $j=1,\ldots,k$ are pairwise disjoint.

Choose a smaller compact neighborhood $E_j\subset\subset  V_j$ 
of $\{x_j,y_j\}$ such that $V_j\bs E_j$ is a union of finitely 
many segments, none of them relatively compact in $V_j$.
(We have three segments for a suitable choice of $E_j$.)
Set $V=\cup_{j=1}^k V_j$ and $E=\cup_{j=1}^k E_j$. 
Then $Q_0 = Q\bs E$ is an open, one dimensional, paracompact,
oriented Hausdorff manifold, hence a union of open segments 
and circles. Any circle in $Q_0$ is a connected
component of $Q$; choosing a nowhere vanishing one-form on it 
(in the correct orientation class) does not affect the  
choices that we shall make on the rest of the set. 
We do the same on any open segment of $Q_0$ which is a
connected component of $Q$. 

It remains to consider those open segments of $Q_0$ which intersect 
at least one of the sets $V_j\bs E_j$. Choose such a segment 
$J$ and an orientation preserving parametrization 
$\phi\colon I=(-1,1)\to J$. Let 
$I'=\{t\in I\colon \phi(t)\in V\}$.
Then $I'$ consists of either one or two subintervals
of $I$, each of them having an end point at $-1$ or at $+1$.
Each of these subintervals  is mapped by $\phi$
onto a segment in one of the sets $V_j\bs E_j$. (The other possibilities 
would require a branch point of $Q_0$, a contradiction.) 

Consider the case when $I'=I_0\cup I_1$  
where $I_0=(-1,a)$, $I_1=(b,1)$ for a pair of 
points $-1< a\le b< 1$. Let $j$ and $l$ be chosen
such that $\phi(I_0)\subset V_j$, $\phi(I_1)\subset V_l$
(we might have $j=l$). Note that $\phi(a)$ is an endpoint 
of $V_j$ and $\phi(b)$ is an endpoint of $V_l$. 
Then $\phi^* \theta_j= d(h_j\circ \phi)$ resp.\ 
$\phi^* \theta_l= d(h_l\circ \phi)$ are nowhere vanishing 
one-forms on $I_0$ resp.\ on $I_1$, both positive 
with respect to the standard orientation of $\R$.
Choose a one-form $\tau$ on $(-1,1)$ which
agrees with the above forms near the respective end points
$-1$ and $+1$ (obviously such $\tau $ exists). 
Then $(\phi^{-1})^*\tau$ is a one-form on $J=\phi(I) \subset Q_0$ 
which agrees with $\theta_j$ in a neighborhood of $E_j$, and 
it agrees with $\theta_l$ in a neighborhood of $E_l$.  
Similarly one deals with the case that $J$ intersects 
only one of the sets $V_j$. Performing this construction for each of the 
finitely many segments $J\subset Q_0$ which intersect $V$ 
we obtain a nowhere vanishing one-form on $Q$.
\end{proof}

\begin{remark}
We wish to elucidate the connection with the methods 
of Straube and Sucheston \cite{SS2}. 
Here we constructed a closed one-form by working directly on 
the leaf space (Proposition \ref{p4.7}).
In \cite{SS2} the authors look instead for a multiplier $h$ 
such that $e^h\eta$ is a closed one-form defining the Levi
foliation $\cL$. To understand this, 
choose a vector field $T$ on $M$ satisfying $\la\eta, T\ra = T\rfloor \eta =1$
and set $\alpha= -{\rm Lie}_T (\eta)$, the  Lie derivative of $\eta$ in the direction $-T$.  
Then $d\eta= \alpha\wedge\eta$, and $\alpha=0$ if $d\eta=0$ \cite[Proposition 2.2]{Ton}. 
The pull-back $\alpha_L$ to a leaf is always closed, and hence 
it defines a De Rham cohomology class $[\alpha_L]\in H^1(L;\R)$. 
From
\[
	d(e^h\eta)= e^h(dh\wedge\eta + d\eta) = e^h(dh+\alpha)\wedge \eta  
\]
we infer that $e^h\eta$ is closed if and only if $dh|_L+\alpha_L=0$ 
on every leaf $L$. A solution on $L$ exists if and only if 
$[\alpha_L]=0\in H^1(L;\R)$, and this holds precisely when
$L$ has trivial infinitesimal holonomy. 
(Indeed, for any closed oriented curve $\gamma\subset L$ we have 
$\oint_\gamma\alpha =\psi_\gamma'(0)-1$, where $t\to \psi_\gamma(t)$ 
is the holonomy map of $\gamma$  on a local transversal to $\cL$ at 
a point $p=\psi_\gamma(0)\in \gamma$.) 
Under suitable geometric conditions on $\cL$, and without assuming that the
foliation extends past $bA$, Straube and Sucheston  \cite{SS2} found a 
solution $h$ when $\dim M=3$.
\end{remark}

%
%
%
%
\section{Levi-flat hypersurfaces without one-sided holonomy} 
The following result shows sufficiency of condition (c) in Proposition \ref{p1.2}.

\begin{theorem}
\label{t5.1}
Let $M$ be a smooth orientable Levi-flat hypersurface
in a complex manifold $X$ such that every element of $\pi_1(M)$ 
is of finite order and its Levi foliation $\cL$ has no nontrivial 
one-sided holonomy (the latter holds in particular if $M$ is real analytic,
or if the leaves of $\cL$ are simply connected).
Then the restriction $\cL|_\omega$ to any open relatively compact subset 
$\omega\subset\subset  M$ is defined by a closed smooth one-form, 
and $M$ admits a $\cC^\infty$ defining function $v$ such that $dd^c v$ 
is flat on $\omega$. If in addition  $\omega$ is simply connected,
there are an open set $U\subset X$, with $M\cap U=\omega$, and a function 
$f=u+iv \in \cC^\infty(U)$ such that $v$ is a defining 
function for $M\cap U$, $u|_\omega$ is a first integral of the 
Levi foliation  $\cL|_\omega$, and the forms $\dibar f$, $dd^c u$ and 
$dd^c v$ are flat on $\omega$.
\end{theorem}

The example in \S 8 shows that Theorem \ref{t5.1} may fail if 
$M$ is real analytic and $\pi_1(M)=\Z$.
In general the function $f$ in Theorem 
\ref{t5.1} cannot be chosen holomorphic even if $M$ 
is real analytic and $\overline\omega$ is contractible \cite{F2}.

\begin{proof}
Let $\cL$ be a smooth, transversely orientable, codimension one
foliation of a connected smooth manifold $M$. By Haefliger \cite{Ha} 
(see also \cite[p.\ 228, Theorem 1.3]{Go1})
the existence of a closed, null homotopic transversal 
$\gamma\subset M$ to $\cL$ implies the existence of a leaf
$L\in\cL$ with nontrivial one-sided holonomy along some
closed curve in $L$. (The proof relies on the Poincar\'e-Bendixson 
theorem   applied to the characteristic foliation induced by 
$\cL$ on a two-disc with boundary $\gamma$.) 
Hence the conditions in Theorem \ref{t5.1} imply that
$\cL$ does not admit any closed transversal.

Let $\ell\subset M$ be a smooth embedded arc transverse to $\cL$. 
We claim that each leaf $L\in\cL$ intersects $\ell$ in at most one point
(and hence $\cL$ is a simple foliation). If not, we find a subarc 
$\tau \subset \ell$ whose  endpoints $p$ and $q$ belong to 
the same leaf $L$. Connecting $q$ to $p$ by an arc 
$\tau'\subset L$ we get a closed loop 
$\lambda =\tau\cdotp \tau'\subset M$. Using the triviality
of the normal bundle to $\cL$ along $\tau'$ we can 
modify $\lambda$ in a small tubular neighborhood
of $\tau'$ into a closed transversal $\widetilde \lambda$ 
to $\cL$ \cite[p.\ 228, 1.2.\ (iii)]{Go1}. Since we have seen 
that a closed transversal does not exist under the stated hypotheses, 
this contradiction proves the claim and shows that $\cL$ is a simple foliation.

Theorem \ref{t4.2} now implies that $\cL$  is defined on any relatively 
compact open subset $\omega\subset \subset M$ by a smooth closed 
one-form $\eta$, and Proposition \ref{p3.1} furnishes a 
smooth defining function $v$ with $dd^c v$ flat on $\omega$.

If $\omega$ is simply connected then $\eta=du$ for 
some $u\in \cC^\infty(\omega)$. Clearly $u$ is constant
on the leaves of $\cL$, i.e., a first integral of $\cL|_\omega$.  
Its asymptotic complexification $f=u+iv$ is a smooth function 
in an open set $U\subset X$ with $U\cap M=\omega$ such that 
$v$ is a defining function for $\omega$ and the forms 
$\dibar f$, $dd^c u$ and $dd^c v$ are flat on $\omega$
(see the proof of Proposition \ref{p3.1}).
\end{proof}

%
%
%
\section{Vanishing of foliation cohomology} 
Let $(M,\cL)$ be a $\cC^r$ foliated manifold. Denote  by 
$\cC^r_\cL$  the sheaf of real $\cC^r$ functions 
on $M$ which are constant on the leaves.
If $\cL$ is the Levi foliation of a Levi-flat hypersurface 
$M$ then $\cC^r_\cL$ is the sheaf of 
real valued CR functions of class $\cC^r$ on $M$. 
The following proves Proposition \ref{p1.2} in case (e).

\begin{theorem}
\label{t6.1}
Let $\cL$ be a transversely orientable codimension one
foliation of class $\cC^{r}$ on a manifold $M$
$(r\in \{2,3,\cdots, \infty, \omega \})$.
If $H^1(M,\cC^{r-1}_\cL)=0$ then $\cL$ is given by a closed, 
nowhere vanishing one-form of class $\cC^{r-1}$.
\end{theorem}

The cohomology group $H^1(M,\cC^{r-1}_\cL)$ can be understood 
either as a \v Cech group, or as a de Rham group, \cite{Ma}.

\begin{proof}
Transverse orientability of $\cL$ implies the existence 
of a $\cC^{r-1}$ vector field $\nu$ which is transverse to $\cL$. 
Choose a transversely oriented $\cC^r$ foliation atlas 
$\{(U_j,\phi_j) \colon j\in J\}$ on $M$ defining $\cL$. 
Write $\phi_j=(\phi'_j, h_j)$ where $h_j$ maps $U_j$ onto 
an open interval $I_j\subset \R$ and $\{h_j=c\}$ are the plaques 
of $\cL|_{U_j}$. For any $i,j\in J$ with $U_{ij} := U_i\cap U_j\ne \emptyset$
we have $h_i= \alpha_{ij}\circ h_j$ on $U_{ij}$ where 
$\alpha_{ij}\colon h_j(U_{ij})\to h_i(U_{ij})$ is a 
$\cC^{r}$ diffeomorphism with positive derivative.
(The collection $\{\alpha_{ij}\}$ is a {\em Haefliger cocycle}
defining $\cL$ \cite{CC}, \cite{Ha}.) We may assume 
$\nu(h_j)>0$ for every $j$. Differentiation gives  
\[ 
	\nu(h_i)= (\alpha'_{ij}\circ h_j)\, \nu(h_j) \quad{\rm on}\ U_{ij}.
\]
This shows that the collection of positive functions 
$b_{ij}=\alpha'_{ij}\circ h_j\in \cC^{r-1}(U_{ij})$ is a one-cocycle 
on the covering $\{U_j\}$ with values in the multiplicative sheaf 
$\cB^{r-1}_\cL$ of positive functions of class $\cC^{r-1}$ which are constant 
on the leaves of $\cL$. The exponential map, 
$\exp \colon \cC^{r-1}_\cL \to \cB^{r-1}_\cL$,
$b \to e^b$, defines an isomorphism between the 
two sheaves (the group operation is additive on the first sheaf
and multiplicative on the second). The hypothesis 
$H^1(M,\cC^{r-1}_\cL)=0$ therefore implies that, after passing 
to a finer $\cL$-atlas, the cocycle $b_{ij}$ is a coboundary, 
$b_{ij}=b_j/b_i$ for some $b_j \in \Gamma(U_j, \cB^{r-1}_\cL)$. 
This gives 
\[ 
	b_i\, \nu(h_i) = b_j\, \nu(h_j) \quad {\rm on}\ U_{ij}.
\]
Since $b_j$ is constant on the plaques $\{h_j=c\} \subset U_j$, 
we have $b_j= \beta_j\circ h_j$ for a unique $\cC^{r-1}$ function 
$\beta_j \colon h_j(U_j) \to\R$. Setting $\alpha_j=\int \beta_j$ 
and $u_j =\alpha_j\circ h_j$ we have 
\[
	\nu(u_j)=(\alpha'_j\circ h_j)\, \nu(h_j) = 
	(\beta_j\circ h_j)\, \nu(h_j) = b_j\, \nu(h_j)>0.
\]
We have thus obtained functions $u_j\in \cC^r_\R(U_j)$ $(j\in J)$ 
which are constant on the plaques in $U_j$ and satisfy 
$\nu(u_j)>0$ on $U_j$ and $\nu(u_i)=\nu(u_j)$ on $U_{ij}$. 
Assuming as we may that the sets $U_{ij}$ are connected,
it follows that the differences $c_{ij}=u_j-u_i$ on $U_{ij}$
are real constants.  Hence the collection of
differentials $du_j$ $(j\in J)$ defines a closed nowhere 
vanishing one-form $\eta$ on $M$ with $\ker\eta=T\cL$. 
\end{proof}

%
%
%
%
\section{Transverse holomorphic vector fields}
In this section we prove Theorems \ref{Neumann}, \ref{app1} and \ref{app2}. 
The key ingredient is the following result on the existence of transverse 
holomorphic vector fields. (Compare with Bedford and Forn\ae ss
\cite[Proposition 6.2]{BeF} and Straube and Sucheston 
\cite[p.\ 156, Proposition]{SS1}.)

\begin{proposition}
\label{p7.1}
Let $M\subset \C^n$ be a Levi-flat hypersurface of class $\cC^r$, $r\ge n+2$,
with a $\cC^r$ defining function $v$ satisfying $j^{r-2}dd^c v|_M=0$.
Let $A$ be a compact strongly pseudoconvex set in $M$ 
(see Theorem \ref{t1.1}). 
There exists a sequence of holomorphic vector fields $\xi_\nu$
in open neighborhoods of $A$ in $\C^n$ $(\nu=1,2,\ldots$) such that
the sequence $\xi_\nu v=\la dv,\xi_\nu\ra$ converges to $1$ 
uniformly on $A$ as $\nu\to\infty$.
\end{proposition}

In particular, the holomorphic vector fields $\xi_\nu$
are transverse to $M$ near $A$.  Recall that a function $v$ with required properties
exists if the Levi foliation of $M$ is defined in a neighborhood of $A$ 
by a nowhere vanishing closed one-form of class $\cC^{r-1}$; see Proposition 1.2
for sufficient conditions.

\begin{proof}
By enlarging $A$ slightly inside $M$ we may assume that Proposition \ref{p3.3} 
applies and hence $A$ is uniformly $\cH$-convex. Choose Stein open sets 
$U_\nu\supset A$ in $\C^n$ $(\nu=1,2,3,\ldots)$ satisfying (\ref{uniform-convex}).
Let $g_j= \di v/ \di  z_j$ for $j=1,\ldots,n$.
Note that ${\di g_j}/{\di \bar z_k} = \di^2 v/\di z_j \di \bar z_k$
is flat to order $r-2$ on $M$; in particular, $g_j$ is a CR function on $M$ near $A$.

\begin{lemma}
\label{l7.2}
There exist holomorphic functions $g_{\nu,j} \in\cO(U_\nu)$
$(j=1, \ldots,n$, $\nu=1,2,\ldots)$ satisfying 
\[
	|g_{\nu,j}(z)-g_j(z)| \le const\cdotp \nu^{-1/2}, \quad z\in A.
\]
\end{lemma}

\begin{proof}
Fix $j\in\{1,\ldots,n\}$ and drop $j$ from the notation.
By the assumption we have $|\dibar g(z)|=o\bigl( d(z,A)^{r-2}\bigr)$
where $d(z,A)=\inf\{|z-w| \colon w\in A\}$. Since every point of
$U_\nu$ has distance at most $c\nu^{-1}$ from $A$ and the volume
of $U_\nu$ is proportional to $\nu^{-1}$ (its thickness in 
the normal direction), we get 
\[
	||\dibar g||_{L^\infty(U_\nu)} = o(\nu^{2-r}),\quad
	||\dibar g||_{L^2(U_\nu)} \le ||\dibar g||_{L^\infty(U_\nu)} \sqrt{Vol(U_\nu)} = o(\nu^{3/2-r}).
\]
By H\"ormander \cite{Ho} there exists a solution $u_\nu$ of
\begin{eqnarray*}
	&& \dibar u_\nu =\dibar g \  \ \hbox{on}\ \ U_\nu, \\
	&& ||u_\nu||_{L^2(U_\nu)} \le const\cdotp ||\dibar g||_{L^2(U_\nu)} = o(\nu^{3/2-r}).
\end{eqnarray*}
Recall 
that every $\cC^1$ function $u$ on $B_\e=\{\zeta \in\C^n \colon |\zeta| <\e\}$ satisfies 
\[
	|u(0)| \le const\cdotp \left( \e^{-n} ||u||_{L^2(B_\e)} +
		\e ||\dibar u||_{L^\infty(B_\e)} \right).
\]
Applying this estimate to  $u_\nu$ at points $z\in A$ with $\e=c\nu^{-1}$ gives 
\[
	|u_\nu(z)| \le \nu^n o(\nu^{3/2-r}) + \nu^{-1}  o(\nu^{2-r}) = o(\nu^{n+3/2-r}).
\]
The function $g_\nu=g-u_\nu$ is then holomorphic on $U_\nu$ for every $\nu\in\N$, 
and assuming $r\ge n+2$ we get the estimate in Lemma \ref{l7.2}.
\end{proof}

Now the proof proceeds as in \cite[p.\ 257]{SS2}.
Let $\cC(A)$ denote the uniform algebra of all continuous complex valued functions
on $A$ endowed with the sup norm, and let $\cA(A)$ be the closure in $\cC(A)$ of the subalgebra 
$\cO(A)|_A$ consisting of the restrictions to $A$ of 
functions holomorphic in open neighborhoods of $A$ in $\C^n$.
Lemma \ref{l7.2} shows that  $g_j \in\cA(A)$ for $j=1,\ldots,n$.
Since these functions have no common zeros on $A$,
a theorem of Rossi \cite[Theorem 2.12]{Rossi} gives 
functions $f_j\in \cA(A)$ satisfying
\[
            1=  \sum_{j=1}^n f_j g_j = \sum_{j=1}^n f_j \frac{\di v}{\di z_j} 
\]
on $A$. Choose sequences of holomorphic functions $f_{\nu,j}$
in small open neighborhoods of $A$ in $\C^n$ such that 
$\lim_{\nu\to\infty} f_{\nu,j}|_A =f_j$ for $j=1,\ldots,n$
(the convergence is uniform on $A$). The sequence of holomorphic vector fields
$
	\xi_\nu =\sum_{j=1}^n f_{\nu,j} \frac{\di}{\di z_j}
$
then satisfies the conclusion of Proposition \ref{p7.1}.
\end{proof}

\smallskip
\textit{Proof of Theorem \ref{Neumann}.}
Let $\Omega \subset \subset \C^n$ be a smoothly bounded 
pseudoconvex domain with a smooth defining function $\rho$,
and let $A\subset b\Omega$ denote the set of infinite type points. 
Boas and Straube proved in \cite[p.\ 227, Theorem]{BS2}
that the Neumann operator $N_q$ and the Bergman projection $P_q$ 
are continuous on the $L_2$-Sobolev space $W^s_{(0,q)}(\Omega)$ 
for $0\le q\le n$ and $s\ge 0$ provided that there is a constant $C>1$,
and for each sufficiently small $\e>0$ a smooth vector field
$\xi_\e$ of type $(1,0)$ in an open set $U_\e \supset A$ in $\C^n$, 
satisfying the following estimates at each point of $A$: 
\begin{itemize}
\item[(1)]  $|{\rm arg} (\xi_\e \rho) | < \e$, $C^{-1}< |\xi_\e \rho| < C$, and 
\item[(2)]  $|\la \di\rho, [\xi_\e, \di/\di\bar z_j]\ra |<\e$ for $j=1,\ldots,n$.
\end{itemize}
(See also \cite{St2005}.)
Assume now that the conditions of Theorem \ref{Neumann} are satisfied.
Let $v$ be an asymptotically pluriharmonic defining function 
for the Levi-flat hypersurface $M$ in a neighborhood of
$A\subset M\cap b\Omega$. Proposition \ref{p7.1} furnishes 
{\em holomorphic} vector fields $\xi_\e$ satisfying condition (1) 
for the function $v$ (instead of $\rho$), and $[\xi_\e, \di/\di\bar z_j]=0$ 
since $\xi_\e$ is holomorphic, so (2) holds. 
Since $A$ is the closure of its relative interior in $b\Omega$, 
the gradients of $v$ and $\rho$ are parallel along $A$,
and we infer that (1) also holds for $\rho$. 
The conclusion now follows from the results of \cite{BS1} and \cite{BS2}.

If $A$ is not the closure of its relative interior then $dv$ and $d\rho$ 
might be completely different on $A$ and we could not infer the condition
(1) for $\rho$. This happens in the worm domain
where $A$ is a single annular leaf.

Note that a desired family of vector fields satisfying (1) and (2) is obtained
in \cite{SS1} and \cite{SS2} under weaker hypotheses on $b\Omega$
near the set of infinite type points. In the case at hand 
the construction in Proposition \ref{p7.1} gives a more precise 
result, although it is based on the same ideas.

\medskip
\noindent {\em Proof of Theorem \ref{app1}.} 
Choose a small open Stein neighborhood $U\subset \C^n$ of $A$ 
and write $U\bs M=U_+\cup U_-$, where $U_{\pm}$ are the 
two sides of $M$ in $U$. We denote by $[M\cap U]^{0,1}$ 
the $(0,1)$-part of the integration current over $M\cap U$.
Let $f$ be a CR function of class $\cC^{k,\alpha}$ in a neighborhood 
of $\overline U\cap M$ in $M$. Then $f[M\cap U]^{0,1}$ is 
a $\dibar$-closed current in $U$; since $U$ is Stein, there exists a 
distribution $F$ on $U$ satisfying $\dibar F=f[M\cap U]^{0,1}$. 
As  $f[M\cap U]^{0,1}$ is supported by $M\cap U$, $F|_{U_\pm}$ is holomorphic 
on $U_\pm$. Since $f$ is of class $\cC^{k,\alpha}$
with $0<\alpha<1$, Theorem 1 in \cite{Ch1} implies that each of 
the functions $F|_{U_\pm}$ extends to a function 
$F_\pm \in \cC^{k,\alpha}(U_{\pm} \cup (M\cap U))$,
and the jump formula $f=F_+ - F_-$ holds on $M\cap U$. 

Let $\xi$ be a transverse holomorphic vector field to $M$
in a neighborhood of $A$, furnished by Proposition \ref{p7.1}.
Choosing a correct orientation we may assume that the flow 
$\phi_t$ of $\xi$ carries a neighborhood of $A$ in $M$
to $U_+$ for small $t>0$, and to $U_-$ for small $t<0$.
For sufficienty small $t>0$ the function 
\[ 
	g_t= F_+\circ \phi_t - F_- \circ\phi_{-t}
\]
is well defined and holomorphic in an open neighborhood $V_t\subset \C^n$ 
of $A$, and $\lim_{t\downarrow 0} g_t|_A = f$ in the $\cC^{k,\alpha}$ topology, 
thus completing the proof. 

If the initial CR function $f$ is of class $\cC^k$ for some integer
$k\ge 1$, the above proof gives approximation by holomorphic 
functions in the $\cC^{k-0}$ sense.

\medskip
\noindent {\em Proof of Theorem \ref{app2}.} 
By Proposition \ref{p7.1} there is a holomorphic vector field
in a neighborhood of $A$ in $\C^n$ which is transverse to $M$.
A Stein neighborhood basis of $\overline{\Omega}$ 
is then furnished by \cite[Lemma 7.3]{BeF} (see also 
\cite[Lemma 1]{FN}), and the Mergelyan approximation property (2) 
follows from  the work of Forn\ae ss and Nagel \cite[Theorem 1]{FN}.

%
%
%
%

\section{A Levi-flat hypersurface with a worm on each side}
The example in this section has been suggested to us 
by J.-E.\ Forn\ae ss; it was used in a related context 
by Bedford and Forn\ae ss \cite{BeF}, and the main idea
can already be seen in the worm domain of Diederich and 
Forn\ae ss \cite{DF}. We shall see that all our results 
fail in this example due to the existence of an annular leaf with 
nontrivial holonomy, thus justifying the hypotheses 
in Theorem \ref{t1.1} and Proposition \ref{p1.2}.
(This example was also discussed in \cite[p.\ 147, Remark 4]{SS1}.)

Denote by $(z,w)$ the coordinates on $\C^*\times \C$.
Let $M \subset \C^*\times \C$ be the real analytic hypersurface 
defined by the following equivalent equations:
\[
   M\colon \quad  \Im (we^{i\log z}) =0 \Longleftrightarrow
    \Im (we^{i\log |z|}) =0.
\]
Indeed, the functions under parentheses differ only 
by the positive multiplicative factor $e^{\arg z}$. The first function 
is multivalued pluriharmonic and hence $M$ 
is Levi-flat. Introducing the 
holomorphic map $\Phi\colon \C^2\to\C^*\times\C$, 
$\Phi(\zeta,t)=\left( e^\zeta, te^{-i\zeta} \right)$,
one sees that $M=\Phi(\C\times \R)$.
The restriction $\Phi(\cdotp,t)$ to a leaf 
$\C\times \{t\} \subset \C\times \R$ 
gives a parametrization of the corresponding leaf 
\[ 
	L_t = \{(e^\zeta, te^{-i \zeta})\colon \zeta\in \C\} \subset M,
	\quad t\in \R
\]
in the Levi foliation $\cL$ of $M$. This parametrization is 
biholomorphic if $t\ne 0$ (so $L_t\simeq \C$) 
while for $t=0$ it is the covering map 
$\C \to L_0=\C^*\times \{0\}$, $\zeta\to (e^\zeta,0)$. 

The line $E=\{(1,s) \colon s\in \R\} \subset M$
is a global transversal for $\cL$, and for every $t\in\R$
we have $L_t \cap E= \{(1,te^{2k\pi})\colon k\in \Z\}$.
The only closed leaf is $L_0=\C^*\times\{0\}$ to which 
all other leaves $L_t$ approach spirally.
Identifying $E$ with $\R$ we see that the holonomy of $L_0$ 
along the positively oriented circle $|z|=1$ is 
$s\to s e^{-2\pi}$. The space of leaves $M/\cL$ 
is the union of a point representing $L_0$ and two closed circles,
each representing the leaves $L_t$ for $t>0$ resp.\ for $t<0$.
Writing $h(z,w)=we^{i\log z}$ we see that the 
closed holomorphic one-form 
\[
	\eta = \frac{d h}{h} = i\frac{dz}{z} + \frac{dw}{w}
\]	
on $\C^*\times \C^*$ defines the Levi foliation of $M\bs L_0$;
there is no such closed one-form in any neighborhood of 
$L_0$ due to nontrivial holonomy.

The hypersurface $M$ divides $\C^*\times\C$ in two connected components
\[
	M_{\pm} =\{(z,w)\in \C^*\times\C \colon  
	\pm \Im (we^{i\log |z|} ) >0\} 
\]
which have the essential properties of a {\em worm domain} 
(Diederich and Forn\-\ae ss \cite{DF}). 
The family of complex annuli
\[
	R_s = \{(z,is)\colon e^{-\pi/2}<|z| < e^{\pi/2} \}, \quad s\in \R
\]	
satisfy the following: 
\begin{itemize}
\item[---] $R_0\subset L_0\subset M$,
\item[---] $bR_s\subset M$ for all $s\in \R$,
\item[---] $R_s\subset M_+$ if $s>0$ and $R_s\subset M_-$ if $s<0$.
\end{itemize}
If $f$ is a holomorphic function in an open neighborhood 
of the annular set
\[
    A_0 =\{(z,w)\in M\colon e^{-\pi/2} \le |z| \le  e^{\pi/2},\ |w|\le 1\} 
\]
then by analytic continuation along the family of annuli
$R_s$, $s\in [-1,1]$, one obtains a holomorphic extension of $f$ 
to a neighborhood of the Levi-flat hypersurface 
$R= \cup_{s\in [-1,1]} R_s$ which therefore belongs to the 
holomorphic hull of $A_0$. Since $R$ intersects both $M_+$ and $M_-$,
we see that Theorem \ref{t1.1} fails. 
Likewise Theorem \ref{t1.4} fails for any 
strongly pseudoconvex domain $D\subset \C^*\times \C$ 
containing $A_0$, and Theorem \ref{app2}
fails for any pseudoconvex domain 
$\Omega\subset\subset \C^*\times \C$ 
with $A_0\subset b\Omega$.

Corollary \ref{c1.6} fails as well which is seen as follows.
With $D$ as above set $D_\pm = D\cap M_\pm$,
$A =\overline{D}\cap M\supset A_0$, $\omega =bD_+\bs A$, 
$K=\overline{D}_-$ and $\Omega = D\bs \widehat K \subset D_+$,
where $\widehat K$ denotes the $\cO(\overline D_+)$-hull of $K$.
The above discussion shows that $R\subset \widehat K$ and 
hence $\Omega$ is a proper subset of $D_+$. 
By  \cite{CL1} and \cite{Lu} every continuous CR function on 
$\omega$ extends holomorphically to $\Omega$; since $\Omega$ is 
pseudoconvex \cite{Sl}, there exists $f\in \cO(\Omega) \cap \cC(\Omega\cup\omega)$
which does not extend holomorphically to $D_+$.

Finally, if $\Omega\subset \subset \C^2$ is a smooth 
pseudoconvex domain with $A_0\subset b\Omega$ and such
that $b\Omega\bs A_0$ is strongly pseudoconvex then
the results of Barrett \cite{Ba3} carry over from the 
standard worm domain and show that the $\dibar$-Neumann 
problem is not globally regular on $\Omega$, in the sense
that the $\dibar$-Neumann operator $N_q$ is not continuous on the
Sobolev spaces $W^s_{0,q}(\Omega)$  \cite[p.\ 147, Remark 4]{SS1}.
For a more precise result see Christ \cite{Christ}.

In the above example the {\em worm phenomen}
is caused by an annulus contained in a leaf
with nontrivial infinitesimal holonomy; there exist
small holomorphic deformations of this annulus with 
boundaries contained in $M$ and such that the interiors 
move to both sides of $M$. We now show that this cannot happen along 
a leaf with trivial (first order) infinitesimal holonomy.
Recall that the holonomy along a loop $\gamma\subset L$ 
is {\em infinitesimally trivial to order $k\in\N$} 
if the associated holonomy map on a local transversal 
$\ell$ to $\cL$ at some point $p\in \gamma$ is of the form 
$x\to x+o(|x|^k)$ in some (and hence any) local coordinate 
$x$ on $\ell$, with $x=0$ corresponding to $p$.

\begin{proposition}
\label{p8.1}
Let $M$ be an orientable Levi-flat hypersurface in a complex manifold 
$X$. Assume that $K\subset X$ is a Stein compact which is contained in 
a leaf $L$ of the Levi foliation of $M$.
If the holonomy of $L$ is infinitesimally trivial to the first order 
along each loop contained in an open neighborhood of $K$ in $L$ then there is a 
holomorphic vector field in a neighborhood of $K$ in $X$ 
which is transverse to $M$.
\end{proposition}

Translating $M$ by the flow of a vector field in 
Proposition \ref{p8.1} gives a Stein neighborhood basis 
of a compact set $A\subset M$ which contains $K$ in its 
relative interior; in particular, 
{\em $M$ has no worm in a neighborhood of $K$}.

\begin{proof}
Replacing $X$ by a small Stein neighborhood of $K$ we
may assume that the leaf $L$ is a closed complex submanifold 
of $X$. Let $\cJ$ denote the sheaf of ideals of $L$.
For a fixed $k\in \N$ let $\cP$ be the subsheaf of $\cJ$
consisting of all germs $h_p\in \cJ_p$ $(p\in L)$ such that 
$dh(p)\ne 0$ and $\Im h|_M = o(|h|^k)$ on a neighborhood of 
$p$ in $M$; for $p\notin L$ we take $\cP_p=\cJ_p = \cO_p$.
Let $\widetilde \cP$ denote the image of $\cP$ in $\cJ/\cJ^{k+1}$.
If the holonomy of $L$ is trivial to order $k$ then the sheaf 
$\widetilde \cP$ admits a section $\tilde h$ 
(Barrett \cite[p.\ 361, Proposition]{Ba1}).  
In our case this holds with $k=1$. Since $X$ is Stein, 
Cartan's Theorem B implies that $\tilde h$ lifts to 
a section $h$ of $\cJ$, i.e., $h$ is a holomorphic 
function vanishing on $L$ to the first order and  
satisfying $\Im h = o(|h|)$ on a neighborhood of 
the leaf $L$ in $M$. Choose a holomorphic vector field 
$\xi$ on $X$ satisfying $\langle dh, \xi \rangle = \sqrt{-1}$ 
on $L$. From $\Im h|_M = o(|h|)$ we see that  
$\xi$ is transverse to $M$ along $L$.
\end{proof}

%
%
%
%
\section{Piecewise real analytic Levi-flat hypersurfaces}
In \S 8 we have seen that the conclusion of 
Theorem \ref{t1.1} may fail in the presence of a leaf with 
nontrivial  infinitesimal holonomy, but 
Proposition \ref{p8.1} gives some hope that leaves with trivial 
infinitesimal holonomy may present no problem.
Hence the following is a reasonable problem.

\begin{problem} 
\label{special-leaves}
Suppose that the Levi foliation of $M$ is simple in the complement of
finitely many leaves $L_1,\ldots,L_k$ and the infinitesimal
holonomy of each $L_j$ is trivial (perhaps to a high order). 
Does it follow that each compact set $A\subset M$ as in 
Theorem \ref{t1.1}  is a Stein compact?
\end{problem}

In this section we give examples indicating that foliations 
of this type are rather common in the class of 
{\em piecewise real analytic Levi-flat hypersurfaces}.
We also show that some of the standard constructions of real codimension 
one foliations, such as the creation of a {\em Reeb component},
{\em turbulization} and {\em spinning}, can be performed in this class.

Let $\Sigma$ be a real analytic manifold. Suppose that 
$\Sigma =D_1\cup\cdots\cup D_m$ where each $D_j$ is a  closed 
domain with real analytic boundary $bD_j$ in $\Sigma$ such that 
for $j\ne k$ the intersection $D_j \cap D_k$ is a union 
of connected components of $bD_j$ and $bD_k$ (possibly empty).
Assume that for each $j$ we are given a real analytic
codimension one foliation $\cF_j$ in an open neighborhood 
of $D_j$ in $\Sigma$ such that $bD_j$ is a union of leaves. 
Let $\cF$ denote the foliation of $\Sigma$ whose restriction 
to $D_j$ equals $\cF_j|_{D_j}$. 
If the foliations on each pair of adjacent domains $D_j$, $D_k$
match up to order $r$ along their common boundary components
then $\cF$ is a {\em transversely piecewise real analytic foliation 
of class $\cC^r$}of $\Sigma$, with real analytic leaves. 
(This is a special case of the {\em tangential gluing of foliations} 
\cite[\S 3.4]{CC}.) Recall that a pair of codimension 
one foliations match up to order $r$ along a common boundary leaf $L$ 
if and only if their holonomies along any closed loop $\gamma\subset L$ 
are tangent to the identity map to order $r$ on the respective sides 
of a local transversal to $L$ at a point in $\gamma$
\cite[p.\ 91, Proposition 3.4.2.]{CC}.  $\cF$ is transversely orientable 
when each of the constituent foliations  is such and their orientations 
match along the boundary leaves.

Suppose now that $\Sigma$ is embedded in a complex manifold $(X,J)$ 
as a real analytic {\em maximal real submanifold}, 
meaning that $T_x X= T_x \Sigma\oplus J (T_x \Sigma)$
for all $x\in \Sigma$. (Every real analytic manifold $\Sigma$ 
admits such an embedding into its complexification. 
If $\dim_\R \Sigma=3$ and $\Sigma$ is orientable then it 
admits a maximal real embedding in $\C^3$ \cite{Au}, \cite{F1}, \cite{Gr}.) 
Let $\cF$ be a real analytic codimension one foliation  of $\Sigma$.
Complexifying a leaf $F\in \cF$ gives a complex submanifold
$\wt F$ in an open neighborhood of $\Sigma$, with 
$\wt F\cap \Sigma=F$ and $\dim_\C \wt F=\dim_\R F$. 
Locally we can perform the complexification uniformly for 
all nearby leaves; restricting our attention to a compact subset $A\subset \Sigma$
we thus find an open neighborhood $U\subset X$ of $A$ and a Levi-flat
real analytic hypersurface $M\subset U$ with Levi foliation 
$\cL = \{\wt F \cap U \colon F\in \cF\}$. Along $\Sigma\cap U \subset M$ 
we have $TM = T\cF \oplus J(T\cF)\oplus N = T\cL\oplus N$
where $N\simeq T\Sigma/T\cF$ is the normal bundle of $\cF$ in $\Sigma$.

Essentially the same construction applies to a piecewise real analytic foliation 
$\cF$ of $\Sigma$ described above. If $\cF$ is transversely $\cC^r$ 
at each of the leaves in $bD_j\cap bD_k$ then the complexifications 
of the individual foliations $\cF_j$ on ${D_j}$ 
match to order $r$ along the common boundary components.
In this way we obtain a Levi foliation $\cL$ on a piecewise real analytic  
$\cC^r$ hypersurface $M=M_1\cup\cdots\cup M_m$, with $\cL|_{M_j}$ 
the complexification of $\cF_j$. The foliated manifolds 
$(M,\cL)$ and $(\Sigma,\cF)$ have the same structure, in particular, 
the same space of leaves.  We now turn to concrete examples.

%
%
\subsection{A Levi-flat hypersurface with a Reeb foliation}
(For the standard $\cC^\infty$ Reeb foliation 
see Examples 1.1.12\ and 3.3.11 in \cite{CC}.) 
Choose integers $n,r\ge 1$. Let $(x,t)$ be coordinates on 
$\R^n\times S^1$. Choose a strictly increasing polynomial function
$\lambda \colon \R\to \R$ satisfying $\lambda(0)=-\pi/2$, 
$\lambda(1)=0$, $\lambda(4)=\pi/2$, whose derivatives up to order 
$r$ vanish at $0,1,4$.  The one-form
\[ 
	\omega= \cos \lambda(\rho) d\rho + \sin \lambda(\rho) dt 
			\qquad (\rho=|x|^2 \in \R_+)
\]
is real analytic and nowhere vanishing on $\R^n\times S^1$. 
From $d\omega \wedge \omega=0$ we infer that $\omega$ determines 
a codimension one foliation $\cF$ of $\R^n\times S^1$.
Its restriction to $T=D^{n}\times S^1$ ($D^n=\{x\in\R^n\colon |x|\le 1\}$)
is a Reeb foliation of $T$ (a {\em Reeb component}) with the only 
closed leaf $F_0= bD^n\times S^1=S^{n-1}\times S^1$ to which all 
other leaves spirally approach. This foliation is not defined by a 
closed one-form in any neighborhood of $F_0$ due to nontrivial holonomy 
along the loops $\{x\}\times S^1$ $(|x|=1)$; this holonomy 
is flat to order $r$. The restriction of $\cF$ to $\{|x|<1\} \times S^1$ 
is a simple foliation given by a fibration over $S^1$.

We decompose the three-sphere $S=S^3$ in a union $T_1\cup T_2$
of two solid tori diffeomorphic to $D^2\times S^1$, with 
$T_1\cap T_2 = bT_1= bT_2 = F_0 \simeq S^1\times S^1$. 
Endowing each $T_j$ with a Reeb foliation $\cF_j$ described 
above one obtains a piecewise real analytic Reeb foliation 
$\cF$ of $S$ whose Reeb components match to order $r\in \N$ 
along the boundary leaf $F_0$. 
(For the $\cC^\infty$ case see \cite[p.\ 93, Example 3.4.4.]{CC}.) 
Embedding $S$ as a real analytic totally real submanifold of $\C^3$ \cite{AR}
and complexifying $\cF$ we obtain a piecewise real analytic 
Levi-flat hypersurface $M\subset\C^3$ whose Levi foliation has the 
structure of the Reeb foliation on $S^3$.
$M$ admits an asymptotically defining function in the complement
of the leaf $L_0=\wt F_0$ (the complexification of the 
torus leaf $F_0 \subset S$), but there is no such function 
near $L_0$ due to nontrivial holonomy.

By Novikov \cite{No} every $\cC^2$ foliation of $S^3$ by surfaces 
contains a Reeb component. According to Barrett \cite{Ba1} the 
Reeb foliation on $S^3$ cannot be realized as the Levi foliation of a smooth 
compact Levi-flat hypersurface, but there are topological realizations
with a corner along the torus leaf.

%
%
\subsection{Turbulization}
Let $\cF$ be the foliation in the previous example, but considered 
now on $\{|x|\le 2\} \times S^1$. 
Since $\lambda(4)=\pi/2$ and the derivatives of $\lambda$ up to order 
$r$ vanish at $4$, the one-form $\omega$ is tangent to $dt$ 
to order $r$ along the torus $T' = \{|x|=2\}\times S^1$,
and hence $\cF$ matches along $T'$ to order $r$ 
with the trivial (horizontal) foliation $\cF_0$ of $\R^n\times S^1$
with leaves $\{t=c\}$. Let $\cF_{turb}$ denote the foliation 
of $\R^n\times S^1$ which equals $\cF$ on $\{|x|\le 2\} \times S^1$ 
and equals $\cF_0$ on $\{|x|\ge 2\} \times S^1$.  This deformation 
of $\cF_0$, known as {\em turbulization} \cite{CC},
can be made in a small tubular neighborhood of any closed transversal 
$\gamma$ in a codimension one foliation and produces 
a new Reeb component along $\gamma$. In the real analytic case 
it can be made by a piecewise real analytic deformation which 
is suitable for complexification, thus giving Levi-flat realizations
of turbulization.  

\subsection{Spinning}
The {\em spinning} modification can be made 
at a boundary component $S\subset bM$ of a foliated 
manifold $(M,\cF)$ provided that every leaf of $\cF$ 
intersects $S$ transversely and the 
induced foliation $\{F\cap S\colon F\in \cF\}$ is 
determined by a closed one-form on $S$ 
\cite[p.\ 84, Example 3.3.B.]{CC}. This modification 
changes $S$ into a closed leaf of a new foliation $\cF_{spin}$ 
which is then suitable for tangential gluing along $S$. 
If all data are real analytic then $\cF_{spin}$ can be made
piecewise real analytic and smooth to a given finite order.
In fact, the Reeb component $\cF$ on $D^n\times S^1$ 
in \S 9.1 is obtained by spinning the trivial foliation with leaves 
$D^n\times \{t\}$ along the boundary $bD^n\times S^1$.

\medskip
{\em Acknowledgements.} 
We thank J.-E.\ Forn\ae ss for telling us about the 
example in \S 8, E.\ Straube for pointing out the connection 
with the $\dibar$-Neumann problem and for very helpful 
discussions, G.\ Tomassini for suggesting Theorem \ref{app1},
and J.\ Mr\v cun, T.\ Ohsawa and M.-C.\ Shaw for 
helpful remarks. We are indebted to the referee
for several suggestions which helped us to improve the exposition.

\bibliographystyle{amsplain}

\end{document}